\documentstyle[11pt]{article}
\title{The Deligne-Simpson problem -- a survey}
\author{Vladimir Petrov Kostov\\ \\ \hspace{7cm}
{\sl To the memory of my mother}} 
\date{}
\newtheorem{tm}{Theorem}
\newtheorem{lm}[tm]{Lemma}

\newtheorem{prop}[tm]{Proposition}
\newtheorem{rem}[tm]{Remark}
\newtheorem{rems}[tm]{Remarks}
\newtheorem{defi}[tm]{Definition}
\newtheorem{ex}[tm]{Example}
\newtheorem{conjecture}[tm]{Open questions}
\newtheorem{conj}[tm]{Conjecture}
\newtheorem{nota}{Notation}
\begin{document}
\maketitle

\begin{abstract}
The Deligne-Simpson problem (DSP) (resp. the weak DSP)   
is formulated like this: {\em give necessary and 
sufficient conditions for the choice of the conjugacy classes 
$C_j\subset GL(n,{\bf C})$ or $c_j\subset gl(n,{\bf C})$ so that there exist 
irreducible (resp. with trivial centralizer) $(p+1)$-tuples  
of matrices $M_j\in C_j$ or $A_j\in c_j$ 
satisfying the equality $M_1\ldots M_{p+1}=I$ or $A_1+\ldots +A_{p+1}=0$}. 
The matrices $M_j$ and $A_j$ are interpreted as monodromy operators of 
regular linear systems and as matrices-residua of Fuchsian ones on Riemann's 
sphere. The present paper offers a survey of the results known up to now 
concerning the DSP.\\  

{\bf Key words:} generic eigenvalues, monodromy operator, 
(weak) Deligne-Simpson problem.

{\bf AMS classification index:} 15A30, 15A24, 20G05
\end{abstract}

\tableofcontents

\section{Introduction}
\subsection{Regular and Fuchsian linear systems on Riemann's sphere}

The problem which is the subject of this paper admits a purely algebraic 
formulation. Yet its importance lies in the analytic theory of systems of 
linear differential equations, this is why we start by considering 
the linear system of ordinary differential equations defined on 
Riemann's sphere:

\begin{equation}\label{system}
{\rm d}X/{\rm d}t=A(t)X
\end{equation}
Here the $n\times n$-matrix $A$ is meromorphic on ${\bf C}P^1$, with poles 
at $a_1$, $\ldots$, $a_{p+1}$; the dependent variables $X$ 
form an $n\times n$-matrix. Without loss of generality we assume that 
$\infty$ is not among the poles $a_j$ and not a pole of the $1$-form 
$A(t)$d$t$. In modern literature the terminology of meromorphic 
connections and sections is often 
preferred to the one of meromorphic linear systems and their solutions 
and there is a 1-1-correspondence between the two languages.

\begin{defi}
The linear system (\ref{system}) is called {\em regular} 
at the pole $a_j$ if its 
solutions have a {\em moderate (or polynomial) 
growth rate} there, i.e. for every 
sector $S$ centered at $a_j$ and of sufficiently small radius and for 
every solution $X$ restricted to the sector 
there exists $N_j\in {\bf R}$ such that $||X(t-a_j)||=O(|t-a_j|^{N_j})$ 
for all $t\in S$.  
System (\ref{system}) is {\em regular} if it is 
regular at all poles $a_j$. 

System (\ref{system}) is {\em Fuchsian} if its poles are logarithmic. 
Every Fuchsian system is regular, see \cite{Wa}.
\end{defi}

\begin{rem}
The opening of the sector $S$ might be $>2\pi$. Restricting to a 
sector is necessary because the solutions are, in general, ramified at the 
poles of the system and by turning around the poles much faster than 
approaching them one can obtain any growth rate.
\end{rem}

A Fuchsian system admits the presentation 

\begin{equation}\label{Fuchs} 
{\rm d}X/{\rm d}t=(\sum _{j=1}^{p+1}A_j/(t-a_j))X~,~
A_j\in gl(n,{\bf C})
\end{equation}
The sum of its {\em matrices-residua} $A_j$ equals 0, i.e. 

\begin{equation}\label{A_j}
A_1+\ldots +A_{p+1}=0
\end{equation}
(recall that there is no pole at $\infty$). 

\begin{rem}
The linear equation (with coefficients meromorphic on ${\bf C}P^1$) 
$\sum _{j=0}^na_j(t)x^{(j)}=0$ is {\em Fuchsian} if $a_j$ has poles of order 
only $\leq n-j$. For linear equations being Fuchsian is equivalent to being 
regular. The best studied Fuchsian equations are the hypergeometric one and 
its generalizations (see \cite{BH}, \cite{L}, \cite{Sa} and \cite{Yo}) and the 
Jordan-Pochhammer equation (see \cite{Ha1}, \cite{Po} and \cite{TaBa}).
\end{rem}

Perform the linear change of the dependent variables 

\begin{equation}\label{W}
X\mapsto W(t)X
\end{equation}
where $W$ is meromorphic on 
${\bf C}P^1$. Most often one requires $W$ to be 
holomorphic and holomorphically invertible for $t\neq a_j$, $j=1,\ldots ,p+1$, 
so that no new singular points appear in the system. As a result of the 
change (\ref{W}) system (\ref{system}) undergoes the 
{\em gauge transformation}:

\begin{equation}\label{gauge}
A\rightarrow -W^{-1}({\rm d}W/{\rm d}t)+W^{-1}AW
\end{equation}
This transformation preserves regularity but, in general, it does not 
preserve being Fuchsian. The only invariant under the group 
of linear transformations (\ref{gauge}) is 
the {\em monodromy group} of the system. 

Set $\Sigma ={\bf C}P^1\backslash \{ a_1,\ldots ,a_{p+1}\}$. 
To define the monodromy group one has to fix a base point 
$a_0\in \Sigma$ and a matrix $B\in GL(n,{\bf C})$. The monodromy group is 
defined only up to conjugacy due to the freedom to choose $a_0$ and $B$.

\begin{defi}
Consider the class 
of homotopy equivalence in $\Sigma$ of a closed contour $\gamma$ with base 
point $a_0$ and bypassing the poles of the system. 
The {\em monodromy operator} of system (\ref{system}) defined by this class 
is the linear operator $M$ acting on the  solution space of the system which 
maps the solution $X$ with $X|_{t=a_0}=B$ into the value of its analytic 
continuation along $\gamma$. Notation: $X\stackrel{\gamma}{\mapsto}XM$. 

The {\em monodromy group} is the subgroup of $GL(n,{\bf C})$ generated 
by all monodromy operators. 
\end{defi}

\begin{rem}\label{anti}
The monodromy group is an {\em antirepresentation} 
$\pi _1(\Sigma )\rightarrow GL(n,{\bf C})$ because one has 

\begin{equation}\label{concat} 
X\stackrel{\gamma _1}{\mapsto}XM_1\stackrel{\gamma _2}{\mapsto}XM_2M_1
\end{equation} 
i.e. the concatenation $\gamma _1\gamma _2$ of the two contours defines the 
monodromy operator $M_2M_1$.   
\end{rem}

One usually chooses a {\em standard set of generators} of $\pi _1(\Sigma )$ 
defined by contours $\gamma _j$, $j=1, \ldots ,p+1$ where  
$\gamma _j$ consists of a segment $[a_0,a_j']$ ($a_j'$ being a point 
close to $a_j$), of a small circumference run counterclockwise (centered 
at $a_j$, passing through $a_j'$ and containing inside no pole of the system 
other than $a_j$), and of the 
segment $[a_j',a_0]$. Thus $\gamma _j$ is freely homotopic to a small loop 
circumventing counterclockwise $a_j$ (and no other pole $a_i$). 
The indices of the poles are 
chosen such that the indices of the contours increase from 1 to $p+1$ 
when one turns around $a_0$ clockwise.  

For the standard choice of the contours the generators $M_j$ satisfy the 
relation 

\begin{equation}\label{M_j}
M_1\ldots M_{p+1}=I
\end{equation}
which can be thought of as the multiplicative analog of (\ref{A_j}) if the 
system is Fuchsian. Indeed, the concatenation of 
contours $\gamma _{p+1}\ldots \gamma _1$ is homotopy equivalent to 0 
and equality (\ref{M_j}) results from $(\ref{concat})$ (see Remark~\ref{anti}).

\begin{rems}\label{MexpA}
1) For a Fuchsian system, if the matrix-residuum $A_j$ has no eigenvalues 
differing by a non-zero 
integer, then the monodromy operator $M_j$ defined as above is conjugate to 
$\exp (2\pi iA_j)$. And it is always true that the eigenvalues $\sigma _{k,j}$ 
of $M_j$ equal $\exp (2\pi i\lambda _{k,j})$ where $\lambda _{k,j}$ are the 
eigenvalues of $A_j$. 

2) If the generators $M_j$ of the monodromy group are defined after a 
standard set of contours $\gamma _j$, then they are conjugate to the 
corresponding 
operators $L_j$ of {\em local monodromy}, i.e. when the poles $a_j$ are 
circumvented counterclockwise along small loops. The operators $L_j$ of 
a regular system can be 
computed (up to conjugacy) algorithmically -- one first makes the system 
Fuchsian at $a_j$ by means of a change (\ref{W}) 
as explained in \cite{Mo} and then carries out the computation 
as explained in \cite{Wa}. Thus $M_j=Q_j^{-1}L_jQ_j$ for some 
$Q_j\in GL(n,{\bf C})$ and the difficulty when computing the monodromy group 
of system (\ref{system}) consists in computing the matrices $Q_j$ which is 
a transcendental problem.
\end{rems} 

\begin{ex}
The Fuchsian system d$X/$d$t=(A/t)X$, $A\in gl(n,{\bf C})$, has two poles -- 
at $0$ and at $\infty$, with matrices-residua equal respectively to $A$ and 
$-A$. Any solution to the system is of the form $X=\exp (A\ln t)G$, 
$G\in GL(n,{\bf C})$. To compute the local monodromy around $0$ one has to 
change the argument of $t$ by $2\pi i$. This results in 
$\ln t\mapsto \ln t+2\pi i$ and $X\mapsto XG^{-1}\exp (2\pi iA)G$, i.e. the 
corresponding monodromy operator equals $G^{-1}\exp (2\pi iA)G$. In the same 
way the monodromy operator at $\infty$ equals $G^{-1}\exp (-2\pi iA)G$.
\end{ex}

\subsection{Formulation of the Deligne-Simpson problem (DSP) and of 
its weak version; generic eigenvalues}

In what follows we write ``tuple'' instead of ``$(p+1)$-tuple''. 
It is natural to state the following realization problem: whether for a given 
tuple of local monodromies (around the poles 
$a_j$) defined up to conjugacy there exists a Fuchsian or at least a regular 
system with such local monodromies. The difficulty is that one must 
have (\ref{M_j}). A similar question can be asked for matrices $A_j$ whose 
sum is 0 (see (\ref{A_j})). The problem can be made more precise:

{\em Give necessary and sufficient conditions on the choice of 
the conjugacy 
classes $C_j\subset GL(n,{\bf C})$ or $c_j\subset gl(n,{\bf C})$ so that 
there exist 
irreducible tuples of matrices $M_j\in C_j$ or $A_j\in c_j$ 
satisfying respectively (\ref{M_j}) or (\ref{A_j}).}

This is the {\em Deligne-Simpson problem (DSP)}. ``Irreducible'' 
means ``with no common proper invariant subspace''. In technical terms 
this means that it is impossible to bring by simultaneous conjugation 
the tuple to a block upper-triangular form 
with the same sizes of the diagonal blocks for all matrices $M_j$ or $A_j$.

\begin{rem}\label{whyirred}
The requirement of irreducibility does not appear in a natural way but 
there are several good reasons for its presence in the 
formulation of the DSP. Firstly, for almost all possible eigenvalues of the 
conjugacy classes the monodromy group is indeed irreducible and knowing that 
the group is irreducible eases the resolution of the problem for such 
eigenvalues (called {\em generic}, see Definition~\ref{nongenrel}).  
Secondly, the answer to the problem for generic eigenvalues depends 
actually not 
on the conjugacy classes but only on the Jordan normal forms which they 
define, see Theorem~\ref{generic}, therefore it is profitable to solve the 
problem first not for general but for generic eigenvalues. 
And thirdly, by restricting oneself  
to the case of generic eigenvalues, one avoids situations like the one 
described in Remark~\ref{undesirable}.
\end{rem}

\begin{rem}\label{multadd}
The name of the problem is motivated by the fact that 
in the multiplicative version (i.e. for matrices $M_j$) 
it was stated in the eighties by P.Deligne (in the 
additive, i.e. for matrices $A_j$, it was stated in the nineties 
by the author) and C.Simpson was the first to obtain 
important results towards its resolution, see \cite{Si1}. The multiplicative 
version is more important because the monodromy group is invariant under 
the action of the group of linear changes (\ref{W}) while the 
matrices-residua of a Fuchsian system are not, see rule (\ref{gauge}) and the 
lines following it. The additive version is technically easier to deal with 
and one can deduce corollaries about the multiplicative one due to 1) of 
Remarks~\ref{MexpA}.
\end{rem}

In what follows we assume that the conjugacy classes $C_j$ (resp. $c_j$) 
satisfy the self-evident condition $\prod \det (C_j)=1$ (resp. 
$\sum$Tr$(c_j)=0$); this condition results from (\ref{M_j}) (resp. from  
(\ref{A_j})). In terms of the eigenvalues $\sigma _{k,j}$  
(resp. $\lambda _{k,j}$) of the matrices from $C_j$ (resp. $c_j$) repeated 
with their multiplicities, these conditions read     

\begin{equation}\label{evs}
\prod _{j=1}^{p+1}\prod _{k=1}^n\sigma _{k,j}=1~~,~~{\rm resp.}~~ 
\sum _{j=1}^{p+1}\sum _{k=1}^n\lambda _{k,j}=0
\end{equation} 

A priori, these are the only conditions that have to be satisfied by the 
eigenvalues of the matrices $M_j$ or $A_j$. 

\begin{defi}\label{nongenrel}
An equality of the form 

\begin{equation}\label{relation}
\prod _{j=1}^{p+1}\prod _{k\in \Phi _j}\sigma _{k,j}=1~~,~~{\rm resp.}~~ 
\sum _{j=1}^{p+1}\sum _{k\in \Phi _j}\lambda _{k,j}=0
\end{equation} 
is called a {\em non-genericity relation};  
the non-empty sets $\Phi _j$ contain one and the same number $<n$ of indices  
for all $j$. Eigenvalues that satisfy none of these relations are called 
{\em generic}. Reducible  
tuples exist only for non-generic eigenvalues (the 
eigenvalues of each diagonal block of a block upper-triangular  
tuple satisfy some non-genericity relation). 
\end{defi}

For non-generic eigenvalues one often encounters situations when there exist 
tuples of matrices $M_j$ or $A_j$ but which are reducible, 
and it is reasonable to give the following definition 
(see also Remark~\ref{whyirred}):

\begin{defi}
The formulation of the 
{\em weak Deligne-Simpson problem} is obtained when in the one 
of the DSP the requirement of irreducibility is replaced by 
the weaker requirement the centralizer of the tuple of matrices 
$A_j$ or $M_j$ to be trivial, i.e. reduced to scalars.
\end{defi}

We say that the DSP (resp. the weak DSP) is {\em solvable}  for a 
given tuple of conjugacy classes $C_j$ or $c_j$ 
if there exists an 
irreducible tuple (resp. a tuple with trivial centralizer) 
of matrices $M_j\in C_j$ satisfying (\ref{M_j}) or of matrices $A_j\in c_j$ 
satisfying (\ref{A_j}). By definition, the DSP is solvable for 
$n=1$. Solvability of the DSP implies automatically the one of the weak DSP.

\begin{rem}\label{undesirable}
If one states the problem of existence of tuples of matrices 
$M_j\in C_j$ or $A_j\in c_j$ satisfying  
respectively condition (\ref{M_j}) or (\ref{A_j}) and with no 
requirement of irreducibility or triviality of the centralizer, then solving 
the problem becomes much harder and the answer to it depends essentially 
on the eigenvalues (not only on the Jordan normal forms). E.g., suppose that  
$p=n=2$ and that two of the matrices $M_j$ (resp. $A_j$) have distinct 
eigenvalues $\sigma _{1,j}\neq \sigma _{2,j}$, $j=1,2$ 
(resp. $\lambda _{1,j}\neq \lambda _{2,j}$) while the third must be scalar 
(i.e. $\sigma _{1,3}=\sigma _{2,3}$, 
resp. $\lambda _{1,3}=\lambda _{2,3}$). Then such triples 
exist exactly if $\sigma _{1,1}\sigma _{1,2}\sigma _{1,3}=1$ or 
$\sigma _{1,1}\sigma _{2,2}\sigma _{1,3}=1$ 
(resp. $\lambda _{1,1}+\lambda _{1,2}+\lambda _{1,3}=0$ 
or $\lambda _{1,1}+\lambda _{2,2}+\lambda _{1,3}=0$). Hence, such triples 
exist exactly if the eigenvalues are not generic. 

A geometric motivation why it is natural to add the condition of triviality 
of the centralizer (if to require irreducibility is too much) 
is given in Subsection~\ref{connectedness}.
\end{rem}

\section{Resolution of the DSP for generic eigenvalues}

\subsection{The quantities $d$ and $r$; Simpson's result\protect\label{SR}}

\begin{defi}\label{JNFdefi}
Call {\em Jordan normal form (JNF) of size $n$} a family 
$J^n=\{ b_{i,l}\}$ ($i\in I_l$, $I_l=\{ 1,\ldots ,s_l\}$, $l\in L$) of 
positive integers $b_{i,l}$ 
whose sum is $n$. Here $L$ is the set of indices of 
eigenvalues (all distinct) and 
$I_l$ is the set of Jordan blocks with the $l$-th eigenvalue, $b_{i,l}$ is the 
size of the $i$-th block with this eigenvalue (for each $l$ fixed 
the sizes of the blocks are listed in decreasing order). E.g. the JNF 
$\{ \{ 2,1\} \{ 4,3,1\} \}$ is of size $11$ and with two eigenvalues to the 
first (resp. second) of 
which there correspond two (resp. three) Jordan blocks, of sizes $2$ and $1$ 
(resp. $4$, $3$ and $1$). An $n\times n$-matrix 
$Y$ has the JNF $J^n$ (notation: $J(Y)=J^n$) if to its distinct 
eigenvalues $\lambda _l$, $l\in L$, there belong Jordan blocks of sizes 
$b_{i,l}$. 
\end{defi}

\begin{nota}
We denote by $J(C)$ (resp. by $J(A)$) 
the JNF defined by the conjugacy class $C$ (resp. the JNF of the matrix $A$). 
By $\{ J_j^n\}$ we denote  
a tuple of JNFs, $j=1$,$\ldots$, $p+1$.
\end{nota}

\begin{nota}
For a conjugacy class $C$ in $GL(n,{\bf C})$ or $gl(n,{\bf C})$ denote by 
$d(C)$ its dimension; recall that it is always even. One has $d(C)=n^2-z(C)$ 
where $z(C)$ is the dimension of the centralizer of a matrix from $C$.

For a matrix $Y$ from $C$ set 
$r(C):=\min _{\lambda \in {\bf C}}{\rm rank}(Y-\lambda I)$. The integer 
$n-r(C)$ is the maximal number of Jordan blocks of $J(Y)$ with one and the 
same eigenvalue. 

Set $d_j:=d(C_j)$ (resp. $d(c_j)$), $r_j:=r(C_j)$ 
(resp. $r(c_j)$). The quantities 
$r(C)$ and $d(C)$ depend only on the JNF $J(Y)=J^n$, not 
on the eigenvalues, so we write sometimes $r(J^n)$ and $d(J^n)$.
\end{nota} 

\begin{rem} 
Recall how to compute $z(C)$ (this is explained in \cite{Ga}). 
It depends only on $J(C)$, not on the eigenvalues of $C$. If 
$J(C)$ is diagonal, with multiplicities of the eigenvalues equal to 
$m_1,\ldots, m_s$, $m_1+\ldots +m_s=n$, 
then one has $z(C)=m_1^2+\ldots +m_s^2$. For a general 
$J(C)=J^n=\{ b_{i,l}\}$ (see Definition~\ref{JNFdefi}) one has 
$z(C)=\sum _l\sum _i(2i-1)b_{i,l}$.
\end{rem}

\begin{prop}\label{d_jr_j}
(C. Simpson, see \cite{Si1}.) The 
following two inequalities are necessary conditions for the solvability 
of the DSP in the case of matrices $M_j$:

\[ \begin{array}{lll}d_1+\ldots +d_{p+1}\geq 2n^2-2&~~~~~&(\alpha _n)\\
{\rm for~all~}j,~r_1+\ldots +\hat{r}_j+\ldots +r_{p+1}\geq n&~~~~~&
(\beta _n)\end{array}\]
\end{prop}

It is shown in \cite{Ko1} that the proposition is true in the case of matrices 
$A_j$ as well. Condition $(\beta _n)$ admits the following generalization 
(see \cite{Ko3} and \cite{Ko4}) which in the case of generic eigenvalues 
coincides with it and which for some non-generic eigenvalues is stronger 
than it:

\begin{prop}\label{betageneralized}
The following inequality is a necessary condition for the solvability of the 
DSP for arbitrary conjugacy classes $C_j$:    

\begin{equation}\label{necbisM} 
\min _{b_j\in {\bf C}^*, b_1\ldots b_{p+1}=1}({\rm rk}(b_1M_1-I)+\ldots 
+{\rm rk}(b_{p+1}M_{p+1}-I))\geq 2n~,~M_j\in C_j
\end{equation}

In the case of conjugacy classes $c_j$ a necessary condition is the 
inequality

\begin{equation}\label{necbisA} 
\min _{b_j\in {\bf C}, b_1+\ldots +b_{p+1}=0}({\rm rk}(A_1-b_1I)+\ldots 
+{\rm rk}(M_{p+1}-b_{p+1}I))\geq 2n~,~A_j\in c_j
\end{equation}
\end{prop}  

The following condition, in general, is not necessary and (as we shall see in 
Subsection~\ref{subsunipnilp}) in most cases it is sufficient 
for the solvability of the DSP, see~\cite{Ko3} and~\cite{Ko5}:

\[ (r_1+\ldots +r_{p+1})\geq 2n~~~~~~~~~~~~~~~~(\omega _n)~~~.\]

The basic result from \cite{Si1} is the following theorem (for the proof 
its author uses the results from \cite{Si2} and \cite{Si3}):
   
\begin{tm}\label{CSimpson}
For generic eigenvalues and when one of the matrices $M_j$ has distinct 
eigenvalues the DSP is solvable for matrices $M_j$ 
if and only if there hold conditions 
$(\alpha _n)$ and $(\beta _n)$.
\end{tm}

It is shown in \cite{Ko1} that the theorem is true in the case of matrices 
$A_j$ as well. Moreover, it remains true both for matrices $M_j$ and $A_j$ 
if one of the matrices has eigenvalues of multiplicity $\leq 2$, not 
necessarily distinct ones, see Theorems 19 and 32 from \cite{Ko2}.

\subsection{Resolution of the DSP for generic eigenvalues and arbitrary JNFs 
of the matrices $M_j$ or $A_j$\protect\label{subsgeneric}}

Theorem \ref{CSimpson} gives, in fact, the necessary and sufficient conditions 
upon $p$ conjugacy classes $C_j$ or $c_j$ so that there exists an 
irreducible $p$-tuple of matrices $M_j\in C_j$ or $A_j\in c_j$. Indeed, for 
almost all such choices of $M_j$ or $A_j$ 
the eigenvalues of the $p+1$ matrices 
(one sets $M_{p+1}=(M_1\ldots M_p)^{-1}$ or $A_{p+1}=-A_1-\ldots -A_p$) 
will be generic and $M_{p+1}$ or $A_{p+1}$ will be with distinct eigenvalues. 

So suppose that there is no condition one of the matrices to be with 
distinct eigenvalues. To formulate the result in this case we need the 
following construction. 
For a given tuple of JNFs $\{ J_j^n\}$ with $n>1$, 
which satisfies conditions 
$(\alpha _n)$ and $(\beta _n)$ and doesn't satisfy condition 
$(\omega _n)$ set $n_1=r_1+\ldots +r_{p+1}-n$. Hence, $n_1<n$ and 
$n-n_1\leq n-r_j$. Define 
the tuple $\{ J_j^{n_1}\}$ as follows: to obtain the JNF 
$J_j^{n_1}$ 
from $J_j^n$ one chooses one of the eigenvalues of $J_j^n$ with 
greatest number $n-r_j$ of Jordan blocks, then decreases  
by 1 the sizes of the $n-n_1$ {\em smallest} Jordan blocks with this 
eigenvalue and deletes the Jordan blocks of size 0. We denote the construction 
by $\Psi :\{ J_j^n\}\mapsto \{ J_j^{n_1}\}$.   

\begin{tm}\label{generic}
For given JNFs $J_j^n$ and for generic 
eigenvalues the DSP is solvable for matrices $A_j$ or $M_j$ 
if and only if the following two conditions hold:

{\em i)} The tuple of JNFs $J_j^n$ satisfies the 
inequality $(\beta _n)$;

{\em ii)} The construction $\Psi$ iterated as long as defined stops 
at a tuple of JNFs $J_j^{n'}$ 
satisfying the inequality $(\omega _{n'})$ or with $n'=1$.
\end{tm}

\begin{rem}
It is true that the result of the theorem does not 
depend on the choice in $\Psi$ of 
an eigenvalue with maximal number of Jordan blocks 
belonging to it, although this is not evident.
\end{rem}

\begin{defi}\label{kappadefi}
N. Katz introduced in \cite{Ka} the quantity 
$\kappa =2n^2-\sum _{j=1}^{p+1}d_j$ and called it {\em index of rigidity} 
of the tuple of conjugacy classes $C_j$ or $c_j$ or of the JNFs 
defined by them.
\end{defi} 

\begin{rem}
If condition $(\alpha _n)$ holds for the JNFs $J_j^n$, then the quantity 
$\kappa$ can take the values $2,0,-2,-4,\ldots$. If there exist irreducible 
tuples of matrices for given conjugacy classes, then there is a 
variety of dimension $2-\kappa$ of two by two non-equivalent representations 
defined by such tuples. In particular, for $\kappa =2$ 
(this case is called {\em rigid}) this variety is of dimension $0$ which means 
that it consists of a finite number of points. It is proved in \cite{Ka} and 
\cite{Si1} that in fact for $\kappa =2$ there is a single irreducible tuple 
defined up to conjugacy. 
\end{rem}  

\begin{lm}\label{kappa1} 
(see \cite{Ko3} and \cite{Ko4}). 
The quantity $\kappa$ is invariant for the 
construction $\Psi$.
\end{lm}  

\begin{lm}\label{kappa2}(see \cite{Ko3}).
For a tuple of JNFs $J_j^n$ satisfying condition $(\omega _n)$ there 
holds condition $(\alpha _n)$ which is a strict inequality and, hence, one has 
$\kappa \leq 0$.
\end{lm}

Lemmas \ref{kappa1} and \ref{kappa2} explain why the necessary condition 
$(\alpha _n)$ does not appear explicitly in the above theorem -- by 
Lemma~\ref{kappa1} it suffices to check that condition $(\alpha _{n'})$ holds 
for the tuple of JNFs $J_j^{n'}$. If inequality $(\omega _{n'})$ 
holds for $\{ J_j^{n'}\}$, then inequality $(\alpha _{n'})$ holds as well and 
is strict (Lemma~\ref{kappa2}). If one has $n'=1$, then 
$(\alpha _{n'})$ holds again, it is an equality and, hence, we are in the 
rigid case. 

The rigid case for matrices $M_j$ has been studied in detail by N. Katz, see 
\cite{Ka}. It is explained there how to construct explicitly irreducible 
tuples of matrices 
$M_j\in C_j$ in the rigid case by means of a middle convolution functor on 
the category of perverse sheaves. An algorithm is given in \cite{Ka} which 
tells whether for given conjugacy classes $C_j$ with $\kappa =2$ (and with 
arbitrary, not necessarily generic eigenvalues) the DSP 
is solvable or not. It is shown in \cite{Ka} that the effect of the 
algorithm upon the JNFs is the same as the one of the construction $\Psi$.

\section{Some particular cases of the DSP}

\subsection{The case of unipotent matrices $M_j$ and of nilpotent 
matrices $A_j$\protect\label{subsunipnilp}}

Suppose that the classes $C_j$ are unipotent and that the classes $c_j$ are 
nilpotent. In this case the DSP and the weak DSP admit an easy formulation. 
The interest in this case is motivated by the fact that the eigenvalues 
are ``the least generic'', i.e. they satisfy all possible 
non-genericity relations. By solving the (weak) DSP for nilpotent or unipotent 
matrices one expects to encounter all possible difficulties that would appear 
in its resolution in the general case. 

\begin{rem}
Condition $(\omega _n)$ (it was introduced in Subsection~\ref{SR}) is 
necessary for the solvability of the DSP when the conjugacy classes $C_j$ 
are unipotent and when the conjugacy classes $c_j$ are nilpotent. Indeed, 
for such conjugacy classes it coincides respectively with conditions 
(\ref{necbisM}) and (\ref{necbisA}) from Proposition~\ref{betageneralized}.
\end{rem}

Define as {\em special} the following cases, when each matrix $A_j$ or $M_j$ 
has Jordan blocks of one and the same size (denoted by $l_j$):

\[ \begin{array}{lllllll}a)&p=3&n=2k&k\in {\bf N},k>1&l_j=2&j=1,2,3,4& \\
b)&p=2&n=3k&k\in {\bf N},k>1&l_j=3&j=1,2,3& \\
c)&p=2&n=4k&k\in {\bf N},k>1&l_1=4&l_2=4&l_3=2\\
d)&p=2&n=6k&k\in {\bf N},k>1&l_1=6&l_2=3&l_3=2\end{array}\]

\begin{rem}
The sizes $(l_1,l_2,l_3)$ from cases b), c) and d) are all positive integer 
solutions to the equation $1/l_1+1/l_2+1/l_3=1$.
\end{rem}

Define as {\em almost special} the four cases obtained from the 
special ones when a couple of blocks of sizes $l_j$ (for the maximal of the 
three or four quantities $l_j$) 
is replaced by a couple of blocks of sizes $l_j+1$ and 
$l_j-1$ while the other blocks remain the same. We list here the sizes of the 
Jordan blocks for the four almost special cases:

\[  \begin{array}{lllll}a')&(3,1,2,\ldots ,2)&(2,\ldots ,2)&(2,\ldots ,2)&
(2,\ldots ,2)\\b')&(4,2,3,\ldots ,3)&(3,\ldots ,3)&(3,\ldots ,3)\\
c')&(5,3,4,\ldots ,4)&(4,\ldots ,4)&(2,\ldots ,2)\\
d')&(7,5,6,\ldots ,6)&(3,\ldots ,3)&(2,\ldots ,2)\end{array}\]

The following theorem can be deduced from Theorem 34 from \cite{Ko5} (for the 
latter's proof the results from \cite{Ko6} were used):

\begin{tm}
1) If condition $(\omega _n)$ holds and if all special and almost 
special cases 
are avoided, then the DSP is solvable for such unipotent or nilpotent 
conjugacy classes.   

2) If condition $(\omega _n)$ holds and if all special cases  
are avoided, then the weak DSP is solvable for such unipotent or nilpotent 
conjugacy classes.
\end{tm}

\begin{rems}
1) For the index of rigidity $\kappa$ (see Definition~\ref{kappadefi}) one has 
$\kappa =0$ in all special cases; these are all cases in which condition 
$(\omega _n)$ holds and one has $\kappa =0$, see \cite{Ko6}, Lemma 3. 

2) If in the definition of the special cases one sets $k=1$, then the DSP is 
solvable.

3) The weak DSP (hence, the DSP as well) is not solvable in the special 
cases. This can be deduced from Theorem 28 from \cite{Ko7}; it is proved 
for matrices $A_j$ in \cite{C-B1} as well. The DSP is not 
solvable in the almost special cases for matrices $A_j$, see \cite{C-B1}. 
\end{rems}

\begin{conj}
The DSP is not solvable for matrices $M_j$ in the almost special cases.
\end{conj}

\subsection{Results concerning the additive version of the DSP}

The additive version of the DSP is completely solved in \cite{C-B1}, for 
arbitrary eigenvalues. It uses the results of the earlier papers \cite{C-B2}, 
\cite{C-B3} and \cite{C-BHo}  
on preprojective algebras and moment map for representations of quivers. 
The answer to the question whether for given conjugacy classes $c_j$ the 
DSP is solvable or not depends on the root system for a Kac-Moody Lie algebra 
with symmetric generalized Cartan matrix constructed after the classes $c_j$. 
Special attention is paid to the rigid case and the DSP is completely solved 
for nilpotent conjugacy classes $c_j$. 
Examples are given which show how the answer  
to the DSP depends (for fixed JNFs $J(c_j)$) 
on the eigenvalues of the classes $c_j$.  

In \cite{Kl} triples of Hermitian matrices $A$, $B$, $A+B$ acting 
on the same $n$-space are considered. 
It is shown there that if $\lambda _i(A)$ are the eigenvalues of $A$ 
listed in the decreasing order, then all relations between the eigenvalues of 
the three matrices (except the trace relation tr$A+$tr$B=$tr$(A+B)$) 
are of the form 

\[ \sum _{k\in K}\lambda _k(A+B)\leq \sum _{i\in I}\lambda _i(A)+
\sum _{j\in J}\lambda _j(B)\] 
where the subsets $I$, $J$ $K$ are precisely those triples for which the 
Schubert cycle $s_K$ is a component in the intersection of Schubert cycles 
$s_I.s_J$. The spectra of the three Hermitian matrices $A$, $B$, $A+B$ 
form a polyhedral convex cone in the space of triple spectra. A recursive 
algorithm is given which generates inequalities describing the cone. 

If $\lambda$, $\mu$, $\nu$ and $V_{\lambda }$, $V_{\mu }$, $V_{\nu }$ denote 
respectively highest weights and the corresponding irreducible 
representations of $GL(V)$, then each tensor product 
$V_{\lambda }\otimes V_{\mu}$ is representable as a sum 
$\sum _{\nu}c_{\lambda \mu}^{\nu}V_{\nu}$ of irreducible representations. 
The results from \cite{Kl} and their refinement from \cite{KnTao} imply that 
the lattice points of the Klyachko cone are precisely the triples of weights 
$\lambda$, $\mu$, $\nu$ with non-zero {\em Littlewood-Richardson coefficient} 
$c_{\lambda \mu}^{\nu}$ (see a survey in \cite{Fu1}). 
The {\em Littlewood-Richardson rule} (see \cite{Fu2}) 
is an algorithm computing these coefficients. A geometric version of this 
rule is the {\em Berenstein-Zelevinsky triangle}, see \cite{BeZe}.

\subsection{The rigid case\protect\label{subsrigid}}

In the case when $\kappa =2$ (the rigid case, see Definition~\ref{kappadefi}) 
irreducible tuples of matrices $M_j\in C_j$ (when they exist) are 
unique up to conjugacy, see \cite{Si1} and \cite{Ka}; from here one easily 
deduces unicity in the additive version of the DSP as well. Such tuples are 
called {\em physically rigid} in \cite{Ka} and {\em linearly rigid} in 
\cite{StVo} and elsewhere. Recall that the 
contribution of N. Katz for the study of the 
rigid case was mentioned at the end of Subsection~\ref{subsgeneric}.  

C. Simpson proves in \cite{Si1} that for $\kappa =2$, if one of the matrices 
$M_j$ has distinct eigenvalues, then one has $p=2$; if the three matrices 
are diagonalizable, then there are only four possibilities for the three 
JNFs. We list them here 
by means of the multiplicity vectors of the eigenvalues of the matrices: 

\[ \begin{array}{llll}
{\rm hypergeometric}&(n-1,1)&(1,\ldots ,1)&(1,\ldots ,1)\\
{\rm odd~family}&(\frac{n+1}{2},\frac{n-1}{2})&
(\frac{n-1}{2},\frac{n-1}{2},1)&(1,\ldots ,1)\\
{\rm even~family}&(\frac{n}{2},\frac{n}{2})&(\frac{n}{2},\frac{n-2}{2},1)&
(1,\ldots ,1)\\{\rm extra~case}&(4,2)&(2,2,2)&(1,1,1,1,1,1)\end{array}\]

O. Gleizer constructs in \cite{Gl} triples of matrices $A_j$ with generic 
eigenvalues and zero sum 
from the above four cases and from another extra case denoted by him by $E_8$. 
The entries of the matrices are ratios of 
products of linear forms in the eigenvalues of the three matrices; the 
non-zero coefficients of these linear forms equal $\pm 1$. The Fuchsian 
systems with three poles on ${\bf C}P^1$ whose matrices-residua are the 
given ones can be considered as the closest relatives of the hypergeometric 
equation of Gauss-Riemann because their triples of spectral flags have 
finitely many orbits for the diagonal action of the general linear group 
in the space of solutions. 

In all four cases scalar products are constructed 
such that the three matrices are self-adjoint w.r.t. them. (For the 
hypergeometric family the results have already been obtained in 
\cite{BH} and Fuchsian systems from the hypergeometric family are equivalent 
to the generalized hypergeometric equations studied in \cite{BH}.) 
In the case when the eigenvalues of the three matrices are real 
this implies that the matrices are real as well. Inequalities upon the 
eigenvalues are given so that these scalar products be positive-definite. 
The inequalities describe non-recursively some faces of the Klyachko cone, see 
the previous subsection.

The scalar products are monodromy invariant complex symmetric bilinear forms 
in the space of solutions of Fuchsian systems with the given matrices as 
residua. The generalized hypergeometric equations have been studied also 
in \cite{O}, in {\em Okubo normal form}. (Okubo shows there that any Fuchsian 
equation can be written in the form $(tI-B){\rm d}X/{\rm d}t=AX$, 
$A,B\in gl(n,{\bf C})$.) For such equations (in Okubo normal 
form) a monodromy invariant Hermitian form has been constructed by 
Y. Haraoka in \cite{Ha3}. These equations have been constructed in \cite{Ha2} 
after the classification of the spectral types of rigid irreducible 
Fuchsian equations has been given in \cite{Y}.

In his paper \cite{Gl} O. Gleizer uses the construction in \cite{MWZ} 
of all indecomposable triple partial flag varieties with finitely many orbits 
for the diagonal action of the general linear group. The construction results 
in a list similar to Simpson's list above, with just one more case, the 
$E_8$ one.  

The existence of over 40 series of rigid triples or quadruples (for 
generic eigenvalues, for both versions of the DSP) is proved in \cite{Ko2}. 
They include all rigid tuples with one 
of the matrices having only eigenvalues of multiplicity $\leq 2$; the last 
condition implies that the tuple consists of $\leq 4$ matrices, see Theorem~22 
from \cite{Ko2}.

\subsection{The multiplicative version of the DSP for unitary matrices}

The DSP for matrices $M_j$, when they are presumed to belong not just 
to $GL(n,{\bf C})$ but to $U(n)$, has been considered in \cite{Bel}, 
\cite{Bi1} and \cite{Bi2}. (In \cite{Bi2} the particular case $n=2$ is 
treated.) In contrast to 
the case of $GL(n,{\bf C})$, when the eigenvalues are generic, and when the 
answer (Theorem~\ref{generic}) is a criterium upon the JNFs and does not 
depend upon the eigenvalues, in the case of 
$U(n)$ the answer depends on the eigenvalues themselves. This answer 
contains a sufficient condition for the existence of a monodromy group 
with given local monodromies; the condition is given 
in terms of non-strict inequalities constructed after the eigenvalues 
and their multiplicity. For each such inequality a condition is given 
whether the validity of the equality is necessary for the existence of the 
monodromy group as well. Similar conditions (in terms of the corresponding 
strict inequalities) are given for the existence of 
irreducible monodromy groups with $M_j\in U(n)$.

In \cite{Bi1} the natural bijective correspondence between the set of 
all equivalence classes of representations 
$\pi _1({\bf C}P^1\backslash \{ a_1,\ldots ,a_{p+1}\} )\rightarrow U(n)$ 
and the set of all isomorphism classes of rank $n$ parabolic stable bundles 
over ${\bf C}P^1$ of parabolic degree zero and $\{ a_1,\ldots ,a_{p+1}\}$  
as the parabolic divisor. The space of equivalence classes of representations 
$\pi _1({\bf C}P^1\backslash \{ a_1,\ldots ,a_{p+1}\} )\rightarrow U(n)$ are 
in one-to-one correspondence with the space of $S$-equivalence classes 
of rank $n$ parabolic semistable bundles of parabolic degree zero, see 
\cite{Si3} and \cite{MeSe}. In these correspondences fixing the conjugacy 
class of the local monodromy around $a_i$ is equivalent to fixing the 
parabolic data at $a_i$. 

In \cite{Bel} another formulation and 
proof of the results of \cite{Bi1} is given as well as 
an algorithm permitting to decide whether a rigid local 
system on ${\bf C}P^1\backslash \{ a_1,\ldots ,a_{p+1}\}$ has finite global 
monodromy; the question has been raised by N. Katz in \cite{Ka}. The methods 
involved in \cite{Bel} (Harder-Narasimhan filtrations) 
are used to strengthen Klyachko's results from \cite{Kl} concerning 
sums of Hermitian matrices. 

It has been observed by Agnihotri and Woodward (see \cite{AW}) and 
independently by Belkale (with the help of Pandharipande) that the DSP for 
unitary matrices is related to quantum cohomology. In \cite{AW} the question 
is raised what the eigenvalues of a product of unitary matrices can be. The 
same question is treated from a symplectic viewpoint in \cite{E}. 

The DSP for an arbitrary compact connected simple simply-connected Lie group 
is considered in \cite{TW}. In most papers cited in this subsection the 
results are related to Gromov-Witten invariants of Grassmanians.  

\subsection{The DSP for finite groups and other results}

In \cite{Vo1} finite groups and their {\em quasi-rigid} generating systems 
are considered. For such a group $G$ the tuple of generators 
$g_i\not\in Z(G)$ with $g_1\ldots g_{p+1}=1$ is called {\em quasi-rigid} 
if for any generators $g_j'$ conjugate respectively to $g_j$ and with 
$g_1'\ldots g_{p+1}'=1$ there exists $g\in G$ such that $g_j'=g^{-1}g_jg$ 
for all $j$. If such generating systems exist, then various criteria permit 
to conclude that certain related groups (e.g. $G/Z(G)$) occur as Galois 
groups over ${\bf Q}(x)$, hence, over ${\bf Q}$ as well. The number $p+1$ 
of generators corresponds to the number of branch points of the geometric 
realizations. The paper shows that the almost simple groups $O^{\pm}_n(2)$, 
and the simple groups $\Omega ^{\pm}_n(2)$ and $Sp_n(2)$ with $n\geq 8$ 
and even, occur as Galois groups over ${\bf Q}$, thus filling a gap left by 
rigidity methods for realizations of groups as Galois groups over ${\bf Q}$ 
for $p+1>3$. 

{\em Belyi triples} are quasi-rigid tuples (introduced by Belyi around 1980, 
see \cite{R} and \cite{Vo1}) which are 
used to realize all classical simple groups as Galois groups over the 
cyclotomic closure of the rationals. For Galois realizations over the 
rationals themselves these triples yield only partial results. In \cite{Vo1} 
the author defines {\em Thompson tuples} as sets of $p+1$ elements 
$\sigma _j\in GL(n,K)$ (where $K$ is a field and $n\geq 3$) such that 
$\sigma _1\ldots \sigma _{p+1}=1$, the group generated by the elements 
$\sigma _j$ being an irreducible subgroup of $GL(n,K)$, and $\sigma _j$ 
being {\em perspectivities} for all $j$, i.e. having eigenspaces of dimension 
$n-1$. Such tuples are studied independently of the groups they generate 
and in comparison with Belyi triples. A criterion is given for Thompson tuples 
to be quasi-rigid. Existence of Thompson tuples with specified characteristic 
polynomials is proved as well as a classification theorem for groups 
generated by Thompson tuples when $K$ is finite and $n>8$. A new construction 
of Belyi triples which achieves partial classification results is introduced. 
It is shown in the paper that $O^{\pm}_n(2)$ and $Sp_n(2)$ for 
$n\geq 6$ and even, have rigid triples of rational generators which implies 
that they are realized as Galois groups over ${\bf Q}$. Belyi triples and 
Thompson tuples in characteristic $0$ 
have been considered by other people as well, 
see Remark~4.1.1 from \cite{DR2}.
  
The following Thompson's conjecture has inspired the papers 
\cite{Vo2},\cite{Vo3} and \cite{DR1}: 
{\em For any fixed finite field ${\bf F}_q$ there exist regular 
Galois realizations 
over the rationals for all but finitely many groups $G({\bf F}_q)$ where 
$G$ is a simple algebraic group of adjoint type defined over ${\bf F}_q$.} 
By the rigidity criterion of Galois theory, this 
yields regular Galois realizations over the rationals for the groups 
$GL_n(q)$ and $U_n(q)$ for $q$ odd, $n=2m+1$ and $m>\varphi (q+1)$ where 
$\varphi$ is Euler's $\varphi$-function. Thus the paper \cite{DR1} settles 
the case of $G=PGL_n$ for $q$ odd while the case of even $n$ 
has been dealt with in \cite{Vo2} and \cite{Vo3} using Thompson tuples. 
Dettweiler and Reiter 
are inspired by Katz's book \cite{Ka}, although they obtain their 
results independently. The results from \cite{DR1} and 
\cite{Vo1} -- \cite{Vo3} can be obtained in an easier way using \cite{DR2} 
and \cite{Vo5}.  

In \cite{DR2} the algorithm of Katz (defined by means of a convolution 
functor, see \cite{Ka}) which tells whether for a given tuple of 
conjugacy classes there exists a rigid tuple of matrices $M_j$, 
is given a purely 
algebraic interpretation. Its analog for the additive version of the DSP 
is also defined. The results are applied to the inverse Galois problem to 
obtain regular Galois realizations over ${\bf Q}$ for families of finite 
orthogonal and symplectic groups: the groups 
$PSO_{2m+1}(q)$, $PGO^+_{2m}(q)$ and $PGO^-_{4m+2}(q)$ and 
in the symplectic case (this is a generalization of a result from 
\cite{ThVo}) the groups $PSp_{2m}(q)$ appear regularly as 
Galois groups over ${\bf Q}$ if $q$ is odd and $m>q$. (For overview of the 
results and for related topics see \cite{MaMat}, \cite{R} and \cite{Vo4}.) 

In the papers cited in this subsection often the results from \cite{F}, 
\cite{FVo} and \cite{Mat} on regular Galois realization of groups 
obtained by rigidity or 
by the braid group action are used. This is also valid for \cite{DebDoE} 
where the DSP for finite groups is also implicitly present, 
see p. 122 there.  

\section{The weak DSP}

The first result we mention in this section is a direct generalization of 
Theorem~\ref{CSimpson}:

\begin{tm}(see \cite{Ko9})
If one of the matrices $A_j$ or $M_j$ is with distinct eigenvalues, then 
conditions $(\alpha _n)$ and $(\beta _n)$ together are necessary and 
sufficient for the solvability of the weak DSP.
\end{tm}

Unlike in the case of the DSP, one cannot allow double eigenvalues in one of 
the matrices (see the lines following Theorem~\ref{CSimpson}) -- 
a triple of nilpotent non-zero $2\times 2$-matrices is 
upper-triangular up to conjugacy, hence, its centralizer is non-trivial 
(it contains each of the matrices). 

For $\kappa =2$ and $\kappa =0$ (see Definition~\ref{kappadefi}) there are 
examples (see the next two subsections) in which conditions {\em i)} and 
{\em ii)} of Theorem~\ref{generic} hold but the weak DSP 
is not solvable for such conjugacy classes. The author was not able to 
find such examples for $\kappa <0$.  

\subsection{An example for $\kappa =0$ (see \cite{Ko7})} 

The following remark will be necessary to 
understand fully the two examples:

\begin{rem}\label{tricky}
In the case of matrices 
$A_j$, if the greatest common divisor $d$ of the 
multiplicities of all 
eigenvalues of all $p+1$ matrices is $>1$, then a non-genericity relation 
results automatically from $\sum$Tr$(c_j)=0$ when the multiplicities of all 
eigenvalues are divided by $d$. In  
the case of matrices $M_j$  
the equality $\prod \sigma _{k,j}=1$ implies that 
if one divides by $d$ the multiplicities of all eigenvalues, then their 
product would equal $\xi =\exp (2\pi ik/d)$, $0\leq k\leq d-1$, 
not necessarily 1, and 
a non-genericity relation holds exactly if $\xi$ is a non-primitive root of 
unity (see Example~\ref{ex1}).
\end{rem}

\begin{ex}\label{ex1} 
Let $p+1=n=4$, let for each $j$ the JNF $J_j^4$ consist 
of two Jordan blocks $2\times 2$, with equal eigenvalues. If the eigenvalues 
of $M_1$, $\ldots$, $M_4$ equal $i$, $1$, $1$ and $1$ 
(resp. $-1$, $1$, $1$, $1$), then 
they are generic (resp. they are non-generic) -- when their multiplicities are 
reduced twice, then their product equals $-1$, a primitive root of $1$ of 
order $2$ (resp. their product equals $1$).
\end{ex}

In order to simplify our example we consider only the case of diagonalizable 
conjugacy classes $C_j$ or $c_j$ (although in \cite{Ko7} the general case is 
treated as well). Suppose that $d>1$ (see Remark~\ref{tricky}) and that 
conditions {\em i)} and {\em ii)} of Theorem~\ref{generic} hold. Suppose that 
the only non-genericity relation if any (satisfied by the eigenvalues of 
the matrices $A_j$) is the one obtained by dividing their 
multiplicities by $d$. In the case of matrices $M_j$ suppose that the only 
non-genericity relation if any is the one obtained by dividing the 
multiplicities by the greatest common divisor of $k$ and $d$, see 
Remark~\ref{tricky}. 

\begin{tm}
1) For such conjugacy classes $c_j$ 
the weak DSP is not solvable for matrices $A_j$.

2) For such conjugacy classes $C_j$ the weak DSP is solvable for matrices 
$M_j$ if and only if $\xi$ is a primitive root of unity.
\end{tm}

\subsection{An example for $\kappa =2$ (see \cite{Ko8})}    

\begin{defi}\label{subordinate}
We say that the conjugacy class $c'$ (in $gl(n,{\bf C})$ 
or in $GL(n,{\bf C})$) is {\em subordinate} to the conjugacy 
class $c$ if $c'$ belongs to the closure of $c$. This means that the  
eigenvalues of $c$ and $c'$ are the same and of the same multiplicities, 
and for each eigenvalue $\lambda _i$ and for each $j\in {\bf N}$ one has 
rk$(A-\lambda _iI)^j\geq$rk$(A'-\lambda _iI)^j$ for $A\in c$, $A'\in c'$. 
If $c'\neq c$, then at least one inequality is strict.
\end{defi}

\begin{defi}\label{good}
A tuple of JNFs is {\em good} if it satisfies conditions i) and ii) 
of Theorem~\ref{generic}.
\end{defi}

\begin{defi}\label{special}
Let $n=ln_1$, $l,n_1\in {\bf N}^*$, $n_1>1$. 
The tuple of conjugacy classes $C_j$ or $c_j$ with $\kappa =2$ 
is called $l$-{\em special} if 
for each class $C_j$ (or $c_j$) there exists a class $C_j'$ (or $c_j'$) 
subordinate to it which is a direct sum of $n_1$ copies of a conjugacy 
class $C_j''\subset GL(l,{\bf C})$ (or $c_j''\subset gl(l,{\bf C})$) 
where the tuple of JNFs $J(C_j'')$ (or $J(c_j'')$) is good and the 
product of the eigenvalues of the classes $C_j''$ equals 1 
(see Remark~\ref{tricky} and Example~\ref{ex1}; for the classes $c_j''$ 
the sum of their eigenvalues is automatically $0$). 
\end{defi}

\begin{rems}\label{c_jrig}
1) The index of rigidity of the tuple of conjugacy classes $c_j''$ 
or $C_j''$ equals 2. Indeed, one has $d(c_j)\geq d(c_j')=(n_1)^2d(c_j'')$ 
and if 
$\sum _{j=1}^{p+1}d(c_j'')\geq 2l^2$, then $\sum _{j=1}^{p+1}d(c_j)\geq 2n^2$, 
i.e. the index of rigidity of the tuple of conjugacy classes $c_j$ 
must be non-positive -- a contradiction. The reasoning holds in the 
case of classes $C_j$ as well. 

2) It follows from the above definition that in the case of matrices $A_j$ the 
eigenvalues of an $l$-special 
tuple of JNFs cannot be generic -- their multiplicities are divisible 
by $n_1$ and, hence, they satisfy a non-genericity relation, 
see Remark~\ref{tricky}. Notice that in the case of matrices $M_j$ 
the divisibility by $n_1$ alone of the multiplicities does not imply that the 
eigenvalues are 
not generic, see Remark~\ref{tricky} and Example~\ref{ex1}. Therefore the 
requirement the product of the eigenvalues of the classes $C_j''$ to equal $1$ 
(see the above definition) is essential.
\end{rems}

\begin{defi}
A tuple of conjugacy classes in $gl(n,{\bf C})$ or $GL(n,{\bf C})$ 
is called {\em special} if it is 
$l$-special for some $l$. If in addition, for this $l$, the classes $c_j''$ or 
$C_j''$ are diagonalizable, then the tuple is called 
{\em special-diagonal}.
\end{defi}

\begin{ex}
For $n>1$ a good tuple of unipotent conjugacy classes in 
$GL(n,{\bf C})$ 
or of nilpotent conjugacy classes in $gl(n,{\bf C})$ is 1-special, hence, 
it is special.
\end{ex}  

\begin{ex}
For $n=9$ the triple of conjugacy classes $c_j$ 
defining the JNFs $\{ \{ 2,2,1,1\} ,\{ 1,1,1\} \}$ for $j=1,2$ and 
$\{ \{ 2,2,1,1\} ,\{ 2,1\} \}$ for $j=3$ is good. Although the multiplicities 
of all eigenvalues are divisible by 3, the triple is not special (a priori 
if it is special, then it is 3-special). Indeed, the JNFs are such that 
the conjugacy classes $c_j''$ from the definition of a special tuple 
must be diagonalizable (for each eigenvalue 
of $c_j$ there are at most two Jordan blocks of size $>1$ and this size is 
actually 2). But then $c_j''$ must have each two eigenvalues, of 
multiplicities 1 and 2, which means that the triple $J(c_1'')$, $J(c_2'')$, 
$J(c_3'')$ is not good.
\end{ex}   

\begin{tm}\label{basicresult}
The weak DSP is not solvable for special-diagonal tuples of 
conjugacy classes.
\end{tm}

It is shown in \cite{Ko8} that conditions {\em i)} and 
{\em ii)} of Theorem~\ref{generic} are necessary 
(for $\kappa =2$) for the solvability of the weak DSP. It seems that the 
above theorem is true if ``special-diagonal'' is replaced by ``special'' 
although the author could not prove it. 

\begin{conjecture}
It would be interesting to know 
whether for $\kappa =2$ special tuples are the only ones for 
which conditions {\em i)} and 
{\em ii)} of Theorem~\ref{generic} hold but the weak DSP is not solvable. 
\end{conjecture}

\section{Geometric motivation of the weak DSP}

\subsection{Analytic deformations of tuples of matrices; correspondence 
between JNFs}

When solving the DSP the author often deforms analytically tuples of 
matrices $A_j$ or $M_j$ with trivial centralizers. Thus its triviality 
occurs as a natural condition and this is the first motivation to include 
it in the formulation of the weak DSP. A more geometric motivation is 
given in the next subsection. 
  
Set $A_j=Q_j^{-1}G_jQ_j$, 
$G_j$ being Jordan matrices. Look for a 
tuple of matrices $\tilde{A}_j$ (whose sum is $0$) of the form 

\[ \tilde{A}_j=(I+\varepsilon X_j(\varepsilon ))^{-1}
Q_j^{-1}(G_j+\varepsilon V_j(\varepsilon ))Q_j
(I+\varepsilon X_j(\varepsilon ))\] 
where 
$\varepsilon \in ({\bf C},0)$ and 
$V_j(\varepsilon )$ are given 
matrices analytic 
in $\varepsilon$; they must satisfy the condition  
tr$(\sum _{j=1}^{p+1}V_j(\varepsilon ))
\equiv 0$; set $N_j=Q_j^{-1}V_jQ_j$. The existence of matrices 
$X_j$ analytic in $\varepsilon$ is deduced from 
the triviality of the centralizer,  
see the proof in \cite{Ko4}. Most often one preserves the 
conjugacy classes of all matrices but one.

If for $\varepsilon \neq 0$ 
small enough the eigenvalues of the matrices $\tilde{A}_j$ are generic, then 
their tuple is irreducible. In a similar way one can deform 
analytically tuples depending on a multi-dimensional parameter.

Given a tuple of matrices $M_j$ with  
trivial centralizer and whose product is $I$, look for matrices 
$\tilde{M}_j$ (whose product is $I$) of the form 

\[ \tilde{M}_j=(I+\varepsilon X_j(\varepsilon ))^{-1}(M_j+
\varepsilon N_j(\varepsilon ))(I+
\varepsilon X_j(\varepsilon ))\] 
where the given matrices $N_j$ depend analytically on 
$\varepsilon \in ({\bf C},0)$ and the product of the determinants of the 
matrices $\tilde{M}_j$ is $1$; one looks for $X_j$ analytic in 
$\varepsilon$. The existence of 
such matrices $X_j$ follows again from the triviality of the centralizer, 
see \cite{Ko4}.

It is often convenient to reduce the resolution of the problem to the case of 
semisimple conjugacy classes. This can be done by choosing appropriate 
matrices $V_j$ or $N_j$ in the above deformations; it would be better to 
speak about {\em specializations} because one has to choose these matrices 
in a special way. Namely, one has to choose them so that the diagonalizable 
matrices $\tilde{A}_j(\varepsilon )$ or 
$\tilde{M}_j(\varepsilon )$ to be (for $\varepsilon \neq 0$) 
from conjugacy classes of least 
possible dimension. It turns out that in this case the JNF $J(A_j)$ or 
$J(M_j)$ changes to its {\em corresponding} diagonal JNF as defined below. 
This correspondence of JNFs has been considered in \cite{Kr}.

\begin{defi}
For a given JNF $J^n=\{ b_{i,l}\}$ define its {\em corresponding} 
diagonal JNF ${J'}^n$. A diagonal JNF is  
a partition of $n$ defined by the multiplicities of the eigenvalues. 
For each $l$ fixed, the set 
$\{ b_{i,l}\}$ is a partition of $\sum _{i\in I_l}b_{i,l}$ and 
$J'$ is the disjoint sum of the dual partitions. 

Two JNFs are said to correspond to one another if they correspond to one 
and the same diagonal JNF.
\end{defi}

\subsection{Connectedness of the moduli spaces of tuples of matrices
\protect\label{connectedness}}

In the present subsection we consider sets of the form 

\[ {\cal V}(c_1,\ldots c_{p+1})=
\{ (A_1,\ldots ,A_{p+1})|A_j\in c_j, A_1+\ldots +A_{p+1}=0\} \]

and    
\[ {\cal W}(C_1,\ldots C_{p+1})=
\{ (M_1,\ldots ,M_{p+1})|M_j\in C_j, M_1\ldots M_{p+1}=I\} \]
or just ${\cal V}$ and ${\cal W}$ for short. 
The following theorem is proved in \cite{Ko9}: 

\begin{tm}\label{basictm}
1) For generic eigenvalues and when one of 
the conjugacy classes $c_j$ (resp. $C_j$) is with distinct eigenvalues the 
set ${\cal V}$ (resp. ${\cal W}$) 
is a smooth and connected variety.

2) For arbitrary eigenvalues, when one of 
the conjugacy classes $c_j$ (resp. $C_j$) is with distinct eigenvalues, and if 
there exist irreducible tuples, then the 
closure of the set ${\cal V}$ (resp. ${\cal W}$) 
is a connected variety. The algebraic closures of these sets coincide with 
their topological closures; these are the closures of the 
subvarieties (the latter are connected) consisting of irreducible 
tuples. 
The singular points of the closures are precisely the tuples 
of matrices with non-trivial centralizers (i.e. not reduced to the scalars). 
\end{tm}

\begin{ex}\label{threestrata}
In the 
case of three diagonalizable non-scalar matrices $A_j$ with $p=2$, 
$n=2$, with eigenvalues respectively $(a,b)$, $(c,d)$, 
$(g,h)$, where the only non-genericity relations are 
$a+c+g=b+d+h=0$, the set ${\cal V}$ 
is a stratified variety with three strata -- $S_0$, $S_1$ and $S_2$. The 
stratum $S_i$ consists of triples which up to conjugacy equal 

\[ A_1=\left( \begin{array}{cc}a&0\\0&b\end{array}\right) ~~,~~
A_2=\left( \begin{array}{cc}c&\varepsilon _i\\ \eta _i&d
\end{array}\right) ~~, ~~
A_3=\left( \begin{array}{cc}g&-\varepsilon _i\\-\eta _i&h\end{array}\right) \]
where $\varepsilon _0=\eta _0=0$, $\varepsilon _1=1$, $\eta _1=0$, 
$\varepsilon _2=0$, $\eta _2=1$. Hence, $S_0$ lies in the closure of $S_1$ 
and $S_2$. There 
are no strata of ${\cal V}$ other than $S_0$, $S_1$ and $S_2$ 
(a theorem by N. Katz (see \cite{Ka}) forbids coexistence of 
irreducible and reducible triples in the rigid case), and ${\cal V}$ 
is connected. However, it is not smooth along $S_0$. Indeed, one can deform 
analytically a triple from $S_0$ into one from $S_1$ and into one from $S_2$; 
the triples from $S_1$ and $S_2$ defining different semi-direct sums, the 
strata $S_1$ and $S_2$ cannot be parts of one and the same smooth variety 
containing $S_0$.  
\end{ex}      

\begin{ex}
Again for $p=n=2$, consider the case when $c_1$ is nilpotent non-scalar and 
$c_2=-c_3$ is with eigenvalues $1$, $2$. Then the  
closure of the variety 
${\cal V}(c_1,c_2,c_3)$ consists of three strata -- $T_0$, $T_1$ and $T_2$. Up 
to conjugacy, the triples from $T_0$ equal  
diag$(0,0)$, diag$(1,-1)$, 
diag$(-1,1)$. The ones from $T_1$ (resp. $T_2$) equal up to conjugacy   

\[ A_1=\left( \begin{array}{cc}0&1\\0&0\end{array}\right) ~~,~~
A_2=\left( \begin{array}{cc}1&-1\\0&2
\end{array}\right) ~~{\rm (resp.}~A_2{\rm )}=\left( \begin{array}{cc}2&-1\\0&1
\end{array}\right) ) ~~,~~A_3=-A_1-A_2~.\]
The stratum $T_0$ lies in the closures of $T_1$ and $T_2$ and like in 
Example~\ref{threestrata}, the closure of 
${\cal V}$ is singular along $T_0$. The variety 
${\cal V}$ itself is not connected because 
$T_0\not\in {\cal V}$ and $T_1$ and $T_2$ are not parts of one and the same 
smooth variety (like $S_1$ and $S_2$ from Example~\ref{threestrata}).
\end{ex}

The reader will find other examples illustrating the stratified structure of 
the varieties ${\cal V}$ or 
${\cal W}$ in \cite{Ko10} and \cite{Ko11}, in particular, 
cases when the dimension of the variety is higher than the expected one 
due to a non-trivial centralizer.

\begin{rems}
1) It would be nice to get rid of the condition one of the classes $c_j$ 
or $C_j$ to be with distinct eigenvalues in parts 1) and 2) of  
Theorem~\ref{basictm} (they are true without this condition in the case 
$\kappa =2$, see \cite{Si1} and \cite{Ka} -- for $\kappa =2$ irreducible 
tuples are unique up to conjugacy and there is no coexistence 
of irreducible and reducible tuples). 

2) It would be interesting to prove the connectedness of the closures 
of the varieties 
${\cal V}$ or ${\cal W}$ without the assumption that there are irreducible 
tuples, see part 2) of Theorem~\ref{basictm}. All examples known to 
the author are of connected closures, see \cite{Ko10} and \cite{Ko11}. 
The connectedness of ${\cal V}$ or ${\cal W}$ implies the one of 
the moduli space of tuples of matrices 
from given conjugacy classes with zero sum or whose product is $I$.

3) Part 2) of Theorem~\ref{basictm} explains the interest in the weak DSP.
\end{rems}

\section{The Riemann-Hilbert problem -- closest relative of the DSP}

The Riemann-Hilbert problem (RHP or Hilbert's 21-st problem) is
formulated as follows: 
 
{\em Prove that for any set of points} $a_1,\ldots ,a_{p+1} \in {\bf C}P^1$
{\em and for any set of matrices} $M_1,\ldots ,M_p \in GL(n, {\bf C})$ {\em
there exists a Fuchsian linear system with poles at} $a_1,\ldots ,a_{p+1} $
{\em for which the corresponding monodromy operators are} $M_1,\ldots ,M_p,
~~M_{p+1} = (M_1\ldots M_p)^{-1} .$ 

Historically, the RHP was first stated for Fuchsian equations, not systems, 
but when one counts the number of parameters necessary to parametrize a 
Fuchsian equation and of the ones necessary to parametrize a monodromy 
group generated by $p$ matrices, one sees that the former, in general, 
is smaller than the latter and one has to allow the presence of additional 
{\em apparent} singularities in the equation, i.e. singularities the 
monodromy around which is trivial. 

The history of the RHP is much longer than 
the one of the DSP and we give as references to it \cite{ArIl}, \cite{AnBo} 
and \cite{Bo1}. We mention here only some of its important events. 

It has been believed for a long time that the RHP has a positive
solution for any $n \in {\bf N}$, after Plemelj in 1908 gives a proof with a 
gap, see \cite{P}.
It nevertheless follows from his proof that if one of the monodromy
operators of system (\ref{system}) is diagonalizable, then system 
(\ref{system}) is equivalent
to a Fuchsian one, see \cite{ArIl}. It also follows that any finitely 
generated subgroup
of $GL(n,{\bf C})$ is the monodromy group of a regular system with prescribed
poles which is Fuchsian at all the poles with the exception, possibly,
of one which can be chosen among them at random.

The RHP has a positive solution for $n=2$ which is
due to Dekkers, see \cite{De}. For $n=3$ the answer is negative, see 
\cite{Bo1}. (This result comes after the gap in Plemelj's proof has been 
detected by A.T. Kohn (see \cite{K}) and Il'yashenko (see \cite{ArIl}).) It
was proved by A.A.Bolibrukh, however, that for $n=3$ the problem has a 
positive answer if
we restrict ourselves to the class of systems with irreducible
monodromy groups, see \cite{Bo1}. Later, the author (see 
\cite{Ko12} and \cite{Ko13}) and 
independently A.A.Bolibrukh (see \cite{Bo2}) proved this result for any $n$.

It is reasonable to reformulate the RHP as follows: 

{\em Find necessary and/or sufficient conditions for the choice of 
the monodromy 
operators $M_1$, $\ldots $, $M_p$ and the points $a_1$, $\ldots $, $a_{p+1}$ 
so that there should exist a Fuchsian system with poles at and only at the 
given points and whose monodromy operators $M_j$ should be  the given ones.} 

Up to now A.A. Bolibruch has found many examples of couples 
(poles, reducible monodromy group) for which the answer to the RHP is 
negative, see \cite{Bo1}, \cite{Bo3}. The RHP discusses the 
question when a given monodromy group 
is realized by a Fuchsian system, therefore it is directly connected with the 
DSP, especially with its additive version. To each of the two problems 
there is an elegant answer formulated in the generic case when the monodromy 
group is irreducible or one of the monodromy matrices has distinct 
eigenvalues. In general, however, the answers to the two problems 
remain essentially different. 

For $n\geq 2$ an 
irreducible monodromy group can be a priori 
realized by infinitely many Fuchsian 
systems, with different tuples of 
conjugacy classes of their matrices-residua. 
In \cite{Ko14} the question is discussed when out of these 
infinitely tuples that are a priori possible, infinitely many are 
in fact not encountered.

Author's address: Universit\'e de Nice -- Sophia Antipolis, Laboratoire 
de Math\'ematiques, Parc Valrose, 06108 Nice Cedex 2, France. 
E-mail: kostov@math.unice.fr   
\end{document}